\documentclass[12pt,twoside]{article}
\usepackage[all]{xy}
\usepackage {amssymb,latexsym,amsthm,amsmath}

\errorcontextlines=0 \numberwithin{equation}{section}
     \newcommand{\tr}{\ensuremath{{}^t\!\!}}
     \newcommand{\C}{\mathfrak C}
     \newcommand{\D}{\mathfrak D}
     \newcommand{\I}{\mathfrak I}
     \newcommand{\R}{\mathfrak R} 
     \newcommand{\Z}{\mathcal Z}
     \newcommand{\Aut}{\textup{Aut}}
     \newcommand{\diag}{\textup{diag}}

       \theoremstyle{plain}
       \newtheorem{theorem}{Theorem}[section]
       \newtheorem{lemma}[theorem]{Lemma}
       
       \newtheorem{proposition}[theorem]{Proposition}

       \newtheorem*{theorem*}{Theorem}
       \theoremstyle{definition}
       
       \numberwithin{equation}{subsection}
       
       \newcommand{\ignore}[1]{}
       
       \newcommand{\mynote}[1]{}

\begin{document}
\setcounter{section}{0}
\setcounter{tocdepth}{2}
\title{\bf Conjugacy Classes of Centralizers in $G_2$ }
\author{Anupam Singh}
\date{}
\maketitle         

\begin{abstract}
Let $G$ be an algebraic group of type $G_2$ over a field $k$ of characteristic $\neq 2,3$.
In this paper we calculate centralizers of semisimple elements in anisotropic $G_2$.
Using these, we show explicitly that there are six conjugacy classes of centralizers in the compact real form of $G_2$.
\end{abstract}

\section{Introduction}
Let $G$ be an algebraic group defined over a field $k$.
Let $X$ be a $G$-space.
Two elements $x,y\in X$ are said to have {\bf same orbit type} if
the isotropy subgroups $G_x$ and $G_y$ are conjugate.
Let $G$ be a compact Lie group acting on a compact manifold $M$.
It was conjectured by Montgomery (ref.~\cite{e}, Problem 45) that there are only finitely many orbit types.
Floyd proved that if $G$ is a torus group acting on a compact orientable manifold then there are only a finite number of distinct isotropy subgroups (ref.~\cite{f}, 4.5).
Using the results of Floyd, Mostow (\cite{m}, Theorem) proved Montgomery's conjecture for a compact Lie group $G$ acting on a compact manifold $M$. 
One can consider the action of a group $G$ on itself by conjugation and ask for the orbit types. 
This case has been considered in~\cite{ku} by R. Kulkarni where the orbit types are called $z$-classes.
In~\cite{s} (Section 3.6 corollary 1 to Theorem 2) Steinberg proved that in a reductive algebraic group $G$ (over algebraically closed field), with char $G$ good, the number of conjugacy classes of centralizers of elements of $G$ is finite.
R. Gouraige studied conjugacy classes of centralizers in $M_n(D)$ (the algebra of endomorphisms of a finite dimensional vector space over a central division algebra) in~\cite{g}.
 
In this paper we consider the question of calculating conjugacy classes of centralizers of elements for anisotropic groups of type $G_2$ over a field of characteristic $\neq 2,3$.
Any group of type $G_2$ over $k$ is precisely the group of $k$-algebra automorphisms of some octonion algebra over $k$ (\cite{se}, Chapter III, Proposition 5, Corollary).   
Anisotropic groups of type $G_2$ over $k$ are given by automorphisms of octonion division algebras over $k$.
We calculate centralizers of elements in $G$ and their conjugacy classes when $G$ is of type $G_2$ and anisotropic over $k$. 
We specifically calculate conjugacy classes of centralizers for compact $G_2$ (anisotropic group $G_2$ over $\mathbb R$) and the number of conjugacy classes of centralizers in this case.
In particular we show this number is finite (Theorem~\ref{maintheorem}).
We plan to continue the study of conjugacy classes of centralizers  in other groups as well.

\section{The Group $G_2$ and some of its Subgroups}

We need description of some subgroups of the group $G_2$ here.
Most of these are described in~\cite{st} section $3$.
Let $k$ be a field and $K$ be an algebraically closed field containing $k$.
We assume that characteristic of $k\neq 2,3$.
Let $\C$ be an octonion algebra over $k$. 
Then the algebraic group $G=\Aut(\C_K)$ where $\C_K=\C\otimes_k K$, is the split, connected, simple algebraic group of type $G_2$. 
Moreover the automorphism group $G$ is defined over $k$.
In fact (\cite{se}, Chapter III, Proposition 5, Corollary), any simple group of type $G_2$ over a field $k$ is isomorphic to the automorphism group of an octonion algebra $\C$ over $k$.
We abuse notation and continue to write the group $G$ of type $G_2$ defined over $k$ as $G=\Aut(\C)$ for some $k$-octonion algebra $\C$.

Let $G=\Aut(\C)$ be an anisotropic group of type $G_2$, where $\C$ is an octonion division algebra over $k$. 
We fix some notation here. 
Let $\D$ be a composition subalgebra of $\C$. 
We denote $G(\C/\D)=\{t\in \Aut(\C)\mid t(x)=x \ \forall x\in \D\}$ and $G(\C,\D)=\{t\in \Aut(\C)\mid t(x)\in\D \ \forall x\in \D\}$, subgroups of $G$.

Let $\D=L$ be a quadratic subalgebra. 
Jacobson described the subgroups $G(\C/L)$ in his paper (ref.~\cite{j}). 
We mention the description of this subgroup here. 
Let $L=k(\gamma)$, where $\gamma^2=c.1 \neq 0$. 
Then $L^{\perp}$ is a left $L$ vector space via the octonion multiplication. 
Also, $h \colon  L^{\perp} \times L^{\perp} \longrightarrow L$ defined by $h(x,y)= N(x,y)+c^{-1}\gamma N(\gamma x,y)$ is a non-degenerate hermitian form on $L^{\perp}$ over $L$. 
Any automorphism $t$ of $\C$ fixing $L$ pointwise induces an $L$-linear map $t|_{L^{\perp}}\colon L^{\perp}\longrightarrow L^{\perp}$. 
Then the group $G(\C/L)$ is isomorphic to the unimodular unitary group 
$SU(L^{\perp},h)$ of the three dimensional space $L^\perp$ over $L$ relative to the hermitian form $h$, via the isomorphism, $t \longmapsto t|_{L^{\perp}}$.
We choose $a,b\in \C$ with $N(a)N(b)\neq 0$ and use Cayley-Dickson doubling process to write $\C=(L\oplus La)\oplus (L\oplus La)b$. 
The nontrivial automorphism of $L$ (denoted as $\alpha\mapsto \bar\alpha$) can be lifted to an automorphism $\rho$ of $\C$ with $\rho^2=1$ where $\rho(\alpha_1 1+\alpha_2 a+\alpha_3 b+\alpha_4 ab)=\bar\alpha_1 1+\bar\alpha_2 a+\bar\alpha_3 b+\bar\alpha_4 ab$.
Then (see~\cite{st}, Proposition 3.5) $G(\C, L)\cong G(\C/L)\rtimes H$ where $H=\langle \rho\rangle$.

Now we describe subgroups $G(\C/\D)$ and $G(\C,\D)$ of $\Aut(\C)$ when $\D=Q$, a quaternion subalgebra. 
We denote the norm $1$ elements of $Q$ by $SL_1(Q)$. 
Then we have, by Cayley-Dickson doubling, $\C=Q\oplus Qa$ for some $a\in Q^{\perp}$ with $N(a)\neq 0$. 
Let $\phi\in \Aut(\C)$ be such that $\phi(x)=x$ for all $x\in Q$. 
Then for $z=x+ya\in\C$, there exist $p\in SL_1(Q)$ such that  $\phi(z)=x+(py)a$ (ref.~\cite{sv}, Proposition 2.2.1). 
We denote such an element of $G(\C/Q)$ by $\R_p$. 
In fact, we have $G(\C/Q)\cong SL_1(Q)$ via the map $\R_p\mapsto p$.
Now let $\psi$ be an automorphism of $Q$.
Define $\widetilde{\psi}\in \Aut(\C)$ by $\widetilde{\psi}(x+ya)=\psi(x)+\psi(y)a$. 
Then one checks easily that $\widetilde{\psi}$ is an automorphism of $\C$ that extends $\psi$ on $Q$. 
These automorphisms form a subgroup of $\Aut(\C)$, which, we will abuse notation and, continue to denote by $\Aut(Q)$. 
We know that any automorphism of a central simple algebra is an inner conjugation by Skolem Noether theorem. 
For any $c\in Q^*$ we have the automorphism $\I_c \colon Q\rightarrow Q$ defined by $x\mapsto cxc^{-1}$. 
Hence $\Aut(Q)\cong PGL_1(Q)\cong Q^*/k^*$ by $\I_c\mapsto c$. 
We then have the subgroup $G(\C,Q)$ of automorphisms of $\C$ leaving $Q$ invariant, $G(\C,Q)\cong G(\C/Q)\rtimes \Aut(Q)$.

Some of these subgroups are conjugate in the group $G$.
We have Skolem-Noether type theorem for composition algebras (\cite{sv}, Corollary 1.7.3) which can be used to describe conjugacy of some of these subgroups.
\begin{theorem}[Skolem Noether type Theorem]\label{skno}
Let $\C$ be a composition algebra and let $\D$ and $\D'$ be subalgebras of the same dimension.
Then, every isomorphism from $\D$ onto $\D'$ can be extended to an automorphism of $C$. 
\end{theorem}
We use this theorem to get following,
\begin{proposition}\label{conjsubgp}
Let $\C$ be an octonion algebra over $k$ and $G=\Aut(\C)$.
Let $\D$ and $\D'$ be two composition subalgebras of $\C$.
Suppose $\D$ and $\D'$ are isomorphic composition subalgebras. 
Then the subgroups $G(\C/\D)$ and $G(\C/\D')$ are conjugate in the group $G$. 
Also the subgroups $G(\C,\D)$ and $G(\C,\D')$ are conjugate in the group $G$.
\end{proposition}
\noindent{\bf Proof :} Let $\tilde\phi$ be the isomorphism of $\D$ to $\D'$.
By Theorem~\ref{skno} $\tilde\phi$ can be extended to an automorphism of $\C$, say $\phi$.
Then, it is easy to check $\phi G(\C/\D)\phi^{-1}=G(\C/\D')$ and $\phi G(\C,\D)\phi^{-1}=G(\C,\D')$.
\qed

Next we calculate centralizers of elements in the group $G$ and it turns out that they are contained in one of the subgroups described above.

\section{Centralizers in Anisotropic $G_2$}

Let $G$ be an anisotropic group of type $G_2$ defined over a field $k$.
Then there exists $\C$, an octonion division algebra over $k$, such that $G(k)\cong \Aut(\C)$.
By abuse of notation we write $G=G(k)=\Aut(\C)$, $k$-points of $G$.
We fix these notation for this section.
We want to calculate centralizer of a semisimple element in $G$. 
If $k$ is perfect (e.g., $char(k)=0$) every element in $G$ is semisimple (\cite{bt}, Corollary 8.5).
Let $t\in G$. 
We denote the centralizer of $t$ in $G$ by $\Z_G(t)=\{g\in G \mid gt=tg\}$.
Since $G$ is anisotropic every element of $G$ leaves a subalgebra of $\C$ fixed pointwise (\cite{st}, Lemma 6.1). 
We denote the subalgebra of fixed points of $t$ by  $\C^t=\{x\in \C \mid t(x)=x\}$.
We note that $t$ restricted to $\C_0$, trace zero space of $\C$, is an element of special orthogonal group of an odd dimensional space (\cite{sv}, Proposition 2.2.2), hence $t$ has a fixed point in $\C_0$ by Cartan-Dieudonne theorem.
This implies that the dimension of $\C^t$ is $\geq 2$.
As the dimension of a composition subalgebra can be $2,4$ or $8$ the subalgebra $\C^t$ is either a quadratic field extension of $k$, a quaternion subalgebra or $t=I$. 
Hence we need to calculate centralizers of elements which are contained in a subgroup $G(\C/L)$ where $L$ is a quadratic field extension or $G(\C/Q)$ where $Q$ is a quaternion subalgebra.
However we have,
\begin{proposition}\label{fixedsub}
Let $G$ be an anisotropic group of type $G_2$ over $k$. 
Let $\C$ be the octonion division algebra over $k$ such that $G=\Aut(\C)$. 
Let $t\in G$. 
Then $\C^t$ is a composition subalgebra of $\C$ and the centralizer $\Z_{G}(t)\subset G(\C,\C^t)$.
\end{proposition}
\noindent {\bf Proof :} The subalgebra $\C^t \subset \C$ is a composition subalgebra as $\C$ is division. 
Let $g\in \Z_{G}(t)$. 
Then, $$
t(g(x))=g(x),~\forall x \in \C^t.
$$ 
Hence $g(x)\in \C^t,~ \forall x\in \C^t$. 
This shows that $g\in G(\C,\C^t)$ and  $\Z_{G}(t)\subset G(\C,\C^t)$.  \qed

\subsection{$\C^t$ is a Quaternion Algebra}

Let $\C^t=Q$ be a quaternion subalgebra of the octonion division algebra $\C$, hence $Q$ itself is division. 
\begin{lemma}
With notation as above, let $t=\R_p \in G(\C/Q)$ for some $p\in SL_1(Q)$ with $p\not\in k$. 
Then $\Z_{G}(t)=\{\R_{p_1}\I_{c_1}\in G(\C,Q) \mid p_1c_1\in L\}$ where $L=k(p)$, a quadratic field extension of $k$.
\end{lemma}
\noindent{\bf Proof :} 
We write $\C=Q\oplus Qb$ for some $b\in Q^{\perp}$. 
Then $t(x+yb)=x+(py)b$. 
From Proposition~\ref{fixedsub} we get, $\Z_{G}(t)\subset G(\C,Q)=\{\R_{p_1}\I_{c_1} \mid p_1\in SL_1(Q),c_1\in Q^*\}$. 
Let $g\in \Z_G(t)$ and let $g=\R_{p_1}\I_{c_1}$.
Then 
$$
gt=tg\Rightarrow \R_{p_1}\I_{c_1}R_p=\R_p\R_{p_1}\I_{c_1}\Rightarrow \R_{p_1c_1pc_1^{-1}}\I_{c_1}=\R_{pp_1}\I_{c_1}
$$ 
and we get, $p_1c_1p=pp_1c_1$, i.e., $p_1c_1\in \Z_{Q}(p)=L$ where $L=k(p)$ is a quadratic field extension of $k$. 
Hence $\Z_G(t)=\{\R_{p_1}\I_{c_1}\mid  p_1c_1\in L, c_1\in Q^*, p_1\in SL_1(Q)\}$.  \qed

Now we consider conjugacy classes of centralizers of these elements in $G$.
\begin{lemma}\label{conjcentraquat}
Let $t,t'\in G=\Aut(\C)$.
Let $t$ and $t'$ leave quaternion subalgebra $Q$ and $Q'$ fixed pointwise, respectively.
Suppose $Q$ and $Q'$ are isomorphic.
Let $t=\R_p$ and $t'=\R_p'$ where $p\in SL_1(Q)$ and $p'\in SL_1(Q')$ and both $p,p'\not\in k$.
Suppose $L=k(p)$ and $L'=k(p')$ are isomorphic field extensions of $k$. 
Then $\Z_G(t)$ is conjugate to $\Z_G(t')$ in $G$.
\end{lemma}
\noindent{\bf Proof : } As $Q$ and $Q'$ are isomorphic, by Theorem~\ref{skno} we have an automorphism $\phi$ of $\C$ such that $\phi|_Q$ is the given isomorphism of $Q$ to $Q'$.
By conjugating the element $t'$ by $\phi$ we may assume $t$ and $t'$ both belong to $G(\C/Q)$.
Let $L$ and $L'$ be isomorphic. 
Then there exists an isomorphism $\psi=\I_c$, conjugation by $c\in Q$, which gives the isomorphism of $L$ to $L'$.
Let $g\in \Z_G(\R_p)$.
From previous lemma $g=\R_{p_1}\I_{c_1}$ with $p_1c_1\in L$. 
Then, 
$$
\psi g\psi^{-1}=\I_c\R_{p_1}\I_{c_1}\I_{c^{-1}}=\R_{cp_{1}c^{-1}}\I_{cc_1c^{-1}}.
$$
We note that $cp_1c^{-1}cc_1c^{-1}=cp_1c_1c^{-1}\in L'$.
This implies $\psi g\psi^{-1}\in\Z_G(t')$ and hence $\Z_G(t)$ is conjugate to $\Z_G(t')$ in $G$. \qed

We note that the center of $SL_1(Q)$ is $\{1,-1\}$ and the element $t$ in $\Aut(\C)$ corresponding to $-1$ is a non-trivial involution (i.e. $t^2=1$).
By a similar calculation as above it is easy to see that $\Z_G(t)=G(\C,Q)$ for $1\neq t$ an involution.
In fact any non-trivial involution in $\Aut(\C)$ correspond to a quaternion subalgebra in this fashion. 
Two involutions are conjugate if and only if the corresponding fixed quaternion subalgebras are isomorphic (\cite{st}, Section 4).
From Proposition~\ref{conjsubgp}, two involutions have their centralizers conjugate if and only if the corresponding fixed quaternion subalgebras are isomorphic. 
We observe that the centralizers corresponding to involutions and other type of elements in Lemma~\ref{conjcentraquat} are not isomorphic hence they can not be conjugate in the group $G$.

\subsection{$\C^t$ is a Quadratic Field Extension}

Now suppose $\C^t=L$, a quadratic field extension of $k$ and $t\in G(\C/L)$. 
We denote the non-trivial automorphism of $L$ by $\alpha\mapsto \bar\alpha$.
Moreover we can find $a,b\in \C$ such that $Q=L\oplus La$ is a quaternion subalgebra and $\C=Q\oplus Qb$. 
Then $L^{\perp}$ is a hermitian space over $L$. 
With respect to the basis $\{a,b,c=ab\}$ we write the subgroup $SU(L^{\perp},h)$ of $G$ as $SU(H)=\{A\in SL(3,L)\mid \tr AH\bar A=H\}$ where $H=\diag\{h(a,a),h(b,b),h(ab,ab)\}$. 
We denote the matrix of $t$ as $A$ with respect to this basis.  
We observe that $t$ leaves a quaternion subalgebra fixed pointwise if and only if $1$ is an eigenvalue  of $A$, i.e., $X-1$ is a factor of the characteristic polynomial $\chi_A(X)$.
Since $\C^t=L$, we may assume that $X-1$ is not a factor of $\chi_A(X)$. 
We have,
\begin{lemma}\label{reg}
With notation as above, let $t\in G(\C/L)$. 
Suppose $1$ is not a root of $\chi_A(X)$.
Then, $\Z_G(t)\subset G(\C/L)$.
\end{lemma}
\noindent{\bf Proof : } Let $t\in G(\C/L)$ be such that $t$ does not fix any point in $L^{\perp}$.  
Let $g\in \Z_G(t)$ then $g\in G(\C,L)\cong SU(L^{\perp},h)\rtimes \langle \rho\rangle$ where $\rho$ is an extension of conjugation on $L$ (\cite{st}, Proposition 3.5). 
We denote the matrix of $t$ in $SU(H)$ by $A$.
Then, either $g\in SU(L^{\perp},h)$ or $g=h\rho$ for some $h\in SU(L^{\perp},h)$. 
Suppose $g=h\rho$.
We denote the matrix of $h$ by $B$ and we denote the action of $\rho$ on elements of $L$ by $\alpha\mapsto \bar\alpha$.
We claim that in this case there is no such element in $\Z_G(t)$. 
We have, 
\begin{eqnarray*}
gt(\alpha_1a+\alpha_2b+\alpha_3c)&=&h\rho t(\alpha_1a+\alpha_2b+\alpha_3c) = h\rho(\alpha_1Aa+\alpha_2Ab+\alpha_3Ac) \\ &=& h(\bar \alpha_1\bar Aa+\bar\alpha_2\bar A b+\bar \alpha_3\bar A c) = \bar \alpha_1B\bar Aa+\bar\alpha_2B\bar Ab+\bar \alpha_3B\bar Ac
\end{eqnarray*}
and 
\begin{eqnarray*}
tg(\alpha_1a+\alpha_2b+\alpha_3c)&=& th\rho (\alpha_1a+\alpha_2b+\alpha_3c)= th(\bar \alpha_1a+\bar \alpha_2 b+\bar \alpha_3c)\\ &=& t(\bar \alpha_1Ba+\bar\alpha_2Bb+\bar \alpha_3 Bc)= \bar \alpha_1ABa+\bar\alpha_2ABb+\bar \alpha_3 ABc.
\end{eqnarray*}
As $g\in \Z_G(t)$, above calculation implies  $AB=B\bar A$ where $B\in SU(H)$. 
In this case the characteristic polynomial $\chi_A(X)=X^3-aX^2+\bar a X-1$ will have  property $\chi_A(X)=\chi_{\bar A}(X)=\overline{\chi_A(X)}$ which implies $a=\bar a$. 
But $\chi_A(X)=X^3-aX^2+a X-1=(X-1)(X^2+(1-a)X+1)$ is reducible and $A$ has $1$ as an eigenvalue. 
Which contradicts the assumption that $1$ is not a root of $\chi_A(X)$. 
Hence $\Z_G(t)\subset G(\C/L)$. \qed

Let $t\in G(\C/L)$ and let $A\in SU(H)$ be the corresponding matrix.
Moreover we assume that $1$ is not a root of $\chi_A(X)$.
Suppose $t$ is semisimple, we have the following cases:
\begin{enumerate}
\item 
\begin{enumerate}
\item The characteristic polynomial of $A$ is irreducible over $L$.
\item The characteristic polynomial has form $(X-\alpha)f(X)$ where $f(X)$ is a quadratic irreducible over $L$.
\item The characteristic polynomial of $A$ has distinct roots over $L$, i.e, $\chi_A(X)=(X-\alpha)(X-\beta)(X-\gamma)$. 
\end{enumerate}
\item The characteristic polynomial of $A$ has two distinct roots over $L$, i.e., $\chi_A(X)=(X-\alpha)(X-\beta)^2$ and minimal polynomial is $(X-\alpha)(X-\beta)$ where $\bar\alpha=\alpha^{-1}$ and $\bar\beta=\beta^{-1}$.
\item The characteristic polynomial of $A$ is $\chi_A(X)=(X-\alpha)^3$ and the minimal polynomial is $(X-\alpha)$ where $\bar\alpha=\alpha^{-1}$.
\end{enumerate}
We calculate centralizers in these cases. 
We observe that whenever $H$ is a connected subgroup of a connected algebraic group $G$ and $t\in H$ such that $\Z_G(t)\subset H$ then $t$ is regular (or strongly regular) in $G$ if and only if it is so in $H$.
Hence $t\in G(\C/L)$ with $\C^t=L$ is regular in $G(\C/L)$ if and only if it is regular in $G$ (Lemma~\ref{reg}). 
This implies that any element of $G(\C/L)$ (which does not have $1$ as an eigenvalue) is a regular element of $G$ in the case $1$, i.e., the connected component of centralizers of such elements are maximal tori.
In the case $3$ the element is scalar $A=\alpha.Id$ and in this case centralizer is whole of $SU(H)$.
We calculate the centralizer in case $2$ below.
\begin{lemma}\label{fixfield}
Let $t\in G=\Aut(\C)$.
With notations as above let $t\in G(\C/L)$ and $A$ be the corresponding element in $SU(H)$.
Suppose $1$ is not an eigenvalue of $A$ and $\chi_A(X)=(X-\alpha)(X-\beta)^2$ and minimal polynomial is $(X-\alpha)(X-\beta)$ where $\bar\alpha=\alpha^{-1}$ and $\bar\beta=\beta^{-1}$.
Then the centralizer of $t$ is isomorphic to a subgroup $\left\{\left (\begin{array}{cc}\det(S)^{-1}&0 \\  0&S \\ \end{array}  \right) \in SU(H)\mid S\in U(W,h|_W) \right\}$ for some $2$-dimensional nondegenerate $L$-subspace $W$ of $L^{\perp}$.
\end{lemma}
\noindent{\bf Proof :} We have $\chi_A(X)=(X-\alpha)(X-\beta)^2$ where $\alpha,\beta\in L$ both not equal to $1$ and the minimal polynomial is $(X-\alpha)(X-\beta)$.  
Then we can choose a basis of $L^{\perp}$ consisting of eigenvectors $v_1,v_2,v_3$ corresponding to eigenvalues $\alpha,\beta,\beta$ of $t$.
We denote the subspace generated by $v_2,v_3$ as $W$.
The matrix of $t$ is diagonal with respect to this basis, $A=\diag\{\alpha,\beta,\beta\}$.
As $1$ is not a root of $\chi_A(X)$ the centralizer $\Z_G(t)\cong \Z_{SU(H)}(A)$ which is isomorphic to
 $\Z_{GL_3(L)}(A)\cap SU(H)=\left\{\left (\begin{array}{cc}\det(S)^{-1}&0 \\  0&S \\ \end{array}  \right) \mid S\in U(W,h|_W) \right\}$. \qed

\section{Conjugacy Classes of Centralizers in Compact $G_2$}

Let $G$ be the compact real group of type $G_2$.
Then (ref.~\cite{b} Chapter V, Section 24.6 Statement (ii) under (c)) $G$ is $\mathbb R$-points of an anisotropic group of type $G_2$. 
We consider the unique anisotropic form of $G_2$ over $\mathbb R$. 
As there is a unique anisotropic Pfister form over $\mathbb R$ of dimension $8$, there is a unique octonion division algebra (up to isomorphism) $\C$ over $\mathbb R$ and $G(\mathbb R)\cong \Aut(\C)$.
We calculate centralizers of elements and there conjugacy classes in this case.
Let $\C$ be the octonion division algebra over $\mathbb R$ and we write, by abuse of notation, $G=G(\mathbb R)=\Aut(\C)$. 
Let $t\in G$.
Let $L$ be a quadratic field extension of $\mathbb R$ left fixed pointwise by $t$ (which is isomorphic to $\mathbb C$).
Then $t\in SU(L^{\perp},h)\cong SU(3)=\{A\in GL_3(\mathbb C)\mid \tr A\bar A=1\}$.
We also note that all quadratic field extensions of $\mathbb R$ contained in $\C$ are isomorphic hence all subgroups of type $SU(L^{\perp},h)$ are conjugate (Proposition~\ref{conjsubgp}) in $G$.
We denote the non-trivial automorphism of $L$ (the complex conjugation in case of $\mathbb C$) by $\alpha\mapsto \bar\alpha$.
We note that every element in this case, $k=\mathbb R$, is semisimple in $G$.
Also every element of $SU(3)$ can be diagonalized in $SU(3)$.
Hence any element $A\in SU(3)$ is of one of the following types:
\begin{itemize}
\item[(a)] The characteristic polynomial of $A$ is $\chi_A(X)=(X-\alpha)(X-\beta)(X-\gamma)$ where $\alpha,\beta$ and $\gamma$ are distinct. 
\item[(b)] The characteristic polynomial of $A$ is $\chi_A(X)=(X-\alpha)(X-\beta)^2$ and the minimal polynomial is $(X-\alpha)(X-\beta)$ where $\alpha\bar \alpha=1=\beta\bar \beta $
\item[(c)] The characteristic polynomial of $A$ is $\chi_A(X)=(X-\alpha)^3$ and the minimal polynomial is $(X-\alpha)$ where $\alpha\bar \alpha=1$. 
\end{itemize}

If $A$ has three distinct roots and none of them is $1$ then by Lemma~\ref{reg} the centralizer of $t$ is contained in $G(\C/L)$.
A simple calculation shows that the centralizer of $A$ in $SU(3)$ is a maximal torus.
Since all maximal tori of $G$ are conjugate (ref.~\cite{k}, Chapter IV, Corollary 4.35) and they are contained in a subgroup isomorphic to $SU(3)$, we get that the element $t$ is strongly regular (those elements of which centralizer is a maximal torus) in $G$.
Moreover the centralizers form one conjugacy class.
If $1$ is an eigenvalue the centralizer may not be contained in $G(\C/L)$.
In fact we have,
\begin{lemma}
With notation as above, let $t\in G(\C/L)\subset G$ where $L\subset \C$ is a quadratic field extension of $\mathbb R$.
In addition, we assume that the characteristic polynomial of $t$ has distinct roots.
Then, $\Z_G(t)\subset G(\C/L)$ if and only if $t$ does not leave any quaternion subalgebra fixed pointwise.
Moreover, if $t$ leaves a quaternion algebra fixed pointwise we have $\Z_G(t)\cong \Z_{G(\C/L)}(t)\rtimes \mathbb Z/2\mathbb Z$. 
\end{lemma}
\noindent{\bf Proof : } If $t\in G(\C/L)$ does not leave any quaternion subalgebra fixed pointwise then the characteristic polynomial of $t$ does not have $1$ as a root.
Hence from Lemma~\ref{reg} it follows that $\Z_G(t)\subset G(\C/L)$.
Now suppose $t\in G(\C/L)$ leaves a quaternion subalgebra $Q$ fixed pointwise.
By using Theorem~\ref{skno} we may assume $Q$ contains $L$.
As $t\in G(\C/L)\cong SU(3)$ the corresponding matrix $A$ can be diagonalized in the subgroup $SU(3)$.
We write the matrix of $t$ as $A=\diag\{1,\alpha,\bar\alpha\}$ for some $\alpha\in L$ with $\alpha\bar\alpha=1$.
Note that $\alpha\neq \bar\alpha$ else $\alpha=\pm 1$ and the roots are not distinct.
We check that $\rho\in G(\C,L)$ conjugates $t$ to $t^{-1}$ (note that $\bar A=A^{-1}$). 
We also find an element $B\in SU(3)$ such that $AB=B\bar A$. 
We take $B=\left (\begin{array}{ccc}-1&0&0 \\  0&0&1 \\ 0&1&0\\\end{array}  \right)$ then one can check that $B\in SU(3)$ and $AB=B\bar A$.
Combined together we get an element $g=h\rho$ such that $g\in Z_G(t)$ and $g^2=1$ where $h\in G(\C/L)$ corresponds to the element $B\in SU(3)$.
We check that $\Z_G(t)\cong \Z_{G(\C/L)}(t)\rtimes \langle g\rangle$.
It is also clear that $\Z_G(t)$ is not contained in $G(\C/L)$. \qed

\noindent If $t$ has $1$ as an eigenvalue the centralizer $\Z_G(t)\cong T\rtimes Z/2\mathbb Z$ where $T$ is a maximal torus in $G$ (which is always contained in a subgroup $SU(3)$).

In the second case with respect to the basis chosen as above $A=\diag\{\alpha,\beta,\beta\}$.
 Suppose neither of $\alpha,\beta$ is $1$ then the centralizer is (ref. Lemma~\ref{fixfield}) 
$$
\Z_G(t)\cong \Z_{SU(3)}(A)\cong \left\{\left (\begin{array}{cc}\det(S)^{-1}&0 \\  0&S \\ \end{array}  \right) \mid S\in U(2) \right\} .
$$
If $1$ is an eigenvalue we have $A=\diag\{1,-1,-1\}$.
In this case $t$ is an involution and the centralizer is contained in $G(\C,L)$.
We note that $\rho$ also centralizes $t$ and hence $\Z_G(t)\cong \Z_{G(\C/L)}(t)\rtimes \langle \rho\rangle\cong U(2)\rtimes \mathbb Z/2\mathbb Z$.

In the third case if $\alpha\neq 1$ centralizer is whole of $SU(3)$.
If $\alpha=1$ the centralizer is whole $G$.

For simplicity we assume $L=\mathbb C$ and make a table of representative elements of each orbit type here.
\vskip3mm
\begin{tabular}{|c|c|}
\hline
{\bf representative element of} & {\bf corresponding centralizer}\\
{\bf an orbit type} & \\ \hline

I & G\\ \hline
$\diag(e^{-i(\theta+\phi)},e^{i\theta},e^{i\phi})$ & maximal torus of $G$ \\ 
$e^{i\theta},e^{i\phi}\neq \pm 1, \pm i$ & \\ \hline
$\diag(1,e^{i\theta},e^{-i\theta})$; $e^{i\theta}\neq \pm 1$ & $T\rtimes \mathbb Z/2\mathbb Z$, $T$ a maximal torus \\  \hline
$\diag(e^{-2i\theta},e^{i\theta},e^{i\theta})$; $e^{i\theta}\neq \pm 1$ & $U(2)$ \\  \hline
$\diag(1,-1,-1)$ & $U(2)\rtimes \mathbb Z/2\mathbb Z$ \\  \hline
$\diag(\omega,\omega,\omega)$; $\omega\neq 1, \omega^3=1$ & $SU(3)$\\ \hline 
\end{tabular}
\vskip3mm

Hence we have proved,
\begin{theorem}\label{maintheorem}
Let $G$ be the anisotropic group of type $G_2$ over $\mathbb R$.
Then there are exactly six orbit types (conjugacy classes of centralizers) as displayed in the table above.
\end{theorem}
\noindent{\bf Acknowledgment :} I thank Maneesh Thakur for his help and several discussions. 
I also thank Ravi Kulkarni for suggesting the problem.


\vskip2mm
\noindent Address: Department of Mathematics\\ The Institute of Mathematical Sciences\\
        C. I. T. Campus, Taramani, Chennai 600113, India\\
        email: anupamk18@gmail.com

\end{document}